\documentclass{article}%
\usepackage{bbm}
\usepackage{epsfig}
\usepackage{graphics,graphicx,amssymb,amsmath,verbatim}
\usepackage{amssymb}
\usepackage{mathrsfs}
\usepackage{amsfonts}
\usepackage{amsmath}
\usepackage{graphicx}
\usepackage{easybmat}
\usepackage{color}%
\providecommand{\U}[1]{\protect \rule{.1in}{.1in}}
\newtheorem{theorem}{Theorem}

\newtheorem{corollary}{Corollary}

\newtheorem{definition}{Definition}

\newtheorem{lemma}{Lemma}

\newtheorem{proposition}{Proposition}
\newtheorem{remark}{Remark}

\allowdisplaybreaks[4]

\begin{document}

\title{On Linear Quadratic Optimal Control of Discrete-Time Complex-Valued Linear Systems}
\author{Bin Zhou\thanks{Center for Control Theory and Guidance Technology, Harbin
Institute of Technology, Harbin, 150001, China. Email: binzhoulee@163.com,
binzhou@hit.edu.cn.}}
\date{}
\maketitle

\begin{abstract}
We study in this paper the linear quadratic optimal control (linear quadratic
regulation, LQR for short) for discrete-time complex-valued linear systems,
which have shown to have several potential applications in control theory.
Firstly, an iterative algorithm was proposed to solve the discrete-time
bimatrix Riccati equation associated with the LQR problem. It is shown that
the proposed algorithm converges to the unique positive definite solution
(bimatrix) to the bimatrix Riccati equation with appropriate initial
conditions. With the help of this iterative algorithm, LQR problem for the
antilinear system, which is a special case of complex-valued linear system,
was carefully examined and three different Riccati equations based approaches
were provided, namely, bimatrix Riccati equation, anti-Riccati equation and
normal Riccati equation. The established approach is then used to solve the
LQR problem for discrete-time time-delay system with one step state delay and
a numerical example was used to illustrate the effectiveness of the proposed methods.

\vspace{0.3cm}


\end{abstract}

\section{Introduction}

Complex-valued linear systems refer to linear systems whose right-hand side
dependent on both the state and its conjugate \cite{zhou17arxiv}. We study
complex-valued linear systems because they have several potential applications
in control theory, for example, describing linear dynamical quantum systems
\cite{zhang17auto} and second-order dynamical systems \cite{zhou17arxiv2}.
Recently, we have studied several analysis and design problems for
complex-valued linear systems, including state response, controllability,
observability, stability, pole assignment, stabilization, linear quadratic
regulation (LQR) and observer design \cite{zhou17arxiv}. We have shown that,
with the help of the so-called bimatrix, results obtained for complex-valued
linear systems are quite analogous to those for normal linear systems
\cite{zhou17arxiv}. Moreover, we have shown that the obtained results include
those for normal linear systems \cite{rugh96book} and antilinear systems
\cite{wdl13aucc,wzls15iet}, which are particular cases of complex-valued
linear systems, as special cases \cite{zhou17arxiv}.

The LQR problem is a fundamental problem in both linear systems theory and
optimal control theory, and has been extensively investigated in the
literature \cite{am07book,dl71tac}. For infinite-time LQR problem, it has been
well known that the solution is completely characterized by the associated
algebraic Riccati equation \cite{am07book,dl71tac,kucera72kyb}. LQR problem
has been extended to several different situations. For example, the LQR problem was solved in \cite{wqls16jfi} for the
so-called antilinear system (which is a special case of the complex-valued
linear systems) and a so-called anti-Riccati equation based solution was established.

With the help of the concept of bimatrix, we have recently solved the LQR
problem for complex-valued linear systems \cite{zhou17arxiv}. It was shown
that the existence of an optimal solution is equivalent to the stabilizability
of the complex-valued linear systems, and is also equivalent to the existence
of positive definite bimatrix to some bimatrix Riccati equation
\cite{zhou17arxiv}. In this paper, based on our early work, we continue to
study the LQR problem for discrete-time complex-valued linear systems. We
first establish an iterative algorithm for solving the discrete-time bimatrix
Riccati equation. The convergence of the algorithm is proven. This iterative
algorithm is not only useful for computing the solution (a bimatrix) to the
bimatrix Riccati equation, but is also helpful in establishing theoretical
results for the sol-called anti-Riccati equation associated with LQR problem
for antilinear systems. Indeed, with such an iterative algorithm, we have
shown that, under the stabilizability assumption, the existence of a solution
to the LQR problem for antilinear systems is equivalent to the existence of a
positive definite solution to the anti-Riccati equation, which closes the gap
in \cite{wqls16jfi} where the existence of a positive definite solution to the
ani-Riccati equation was not guaranteed. We will also establish another normal
Riccati equation based solution to the LQR problem for antilinear systems. The
relationships among the bimatrix Riccati equation, anti-Riccati equation and
normal Riccati equation are revealed. At the same time, we show that the
anti-Riccati equation can be equivalently transformed into a nonlinear matrix
equation that has been carefully studied in our early work
\cite{lzl14amc,zcl13amc}. Finally, by expressing a discrete-time linear
time-delay system as a complex-valued system model, the LQR problem for such a
system is solved by using bimatrix Riccati equations. A numerical example was
worked out to illustrate the effectiveness of the proposed approach.

\textbf{Notation}: For a matrix $A\in \mathbf{C}^{n\times m},$ we use $A^{\#},$
$A^{\mathrm{T}},$ $A^{\mathrm{H}},$ $\mathrm{rank}\left(  A\right)  ,$
$\left \Vert A\right \Vert ,$ $\operatorname{Re}\left(  A\right)  $ and
$\operatorname{Im}\left(  A\right)  $ to denote respectively its conjugate,
transpose, conjugate transpose, rank, norm, real part and imaginary part. Thus
$A^{-\#}$ denotes $(A^{\#})^{-1}$ or $(A^{-1})^{\#}.$ Denote $\mathrm{j}$ the
unitary imaginary number. For a matrix pair $\left(  A_{1},A_{2}\right)
\in \left(  \mathbf{C}^{n\times m},\mathbf{C}^{n\times m}\right)  ,$ the
bimatrix $\left \{  A_{1},A_{2}\right \}  \ $is defined in such a manner that
$\left \{  A_{1},A_{2}\right \}  x=A_{1}x+A_{2}^{\#}x^{\#}.$ Further definitions
and properties about bimatrix can be found in \cite{zhou17arxiv}.

\section{\label{sec2}Optimal Control of Complex-Valued Linear Systems}

\subsection{A Brief Introduction to Complex-Valued Linear Systems}

We continue to study in this paper the following complex-valued linear system
\cite{zhou17arxiv,zhou17arxiv2}%
\begin{equation}
x\left(  k+1\right)  =\left \{  A_{1},A_{2}\right \}  x(k)+\left \{  B_{1}%
,B_{2}\right \}  u(k), \label{sys}%
\end{equation}
where $A_{i}\in \mathbf{C}^{n\times n}$ and $B_{i}\in \mathbf{C}^{n\times
m},i=1,2,$ are known coefficients, $x(k)$ is the state, and $u(k)$ is the
control. The initial condition is set to be $x\left(  0\right)  =x_{0}%
\in \mathbf{C}^{n}.$ Clearly, system (\ref{sys}) becomes the normal linear
system
\begin{equation}
x\left(  k+1\right)  =A_{1}x(k)+B_{1}u(k), \label{normal}%
\end{equation}
if $A_{2}$ and $B_{2}$ are null, and becomes the so-called antilinear system
\begin{equation}
x\left(  k+1\right)  =A_{2}^{\#}x^{\#}(k)+B_{2}^{\#}u^{\#}(k),
\label{antilinear}%
\end{equation}
if $A_{1}$ and $B_{1}$ are zeros. The antilinear system (\ref{antilinear}) was
firstly studied in \cite{wdl13aucc} and \cite{wzls15iet}. We have shown
recently in \cite{zhou17arxiv} and \cite{zhou17arxiv2} that the complex-valued
linear system has several potential applications in control, for example, for
control of linear dynamical quantum systems \cite{zhang17auto} and
second-order dynamical systems \cite{zhou17arxiv2}.

In this paper, based on our early work \cite{zhou17arxiv}, we continue to
study the linear quadratic optimal control problem for system (\ref{sys}). To
this end, we introduce some based concepts for this system.

\begin{definition}
\cite{zhou17arxiv} The complex-valued linear system (\ref{sys}) is said to be
stabilizable if there exists a so-called full state feedback
\begin{equation}
u\left(  k\right)  =\left \{  K_{1},K_{2}\right \}  x\left(  k\right)
=K_{1}x\left(  k\right)  +K_{2}^{\#}x^{\#}\left(  k\right)  ,
\label{eqfeedback}%
\end{equation}
such that the following closed-loop system is asymptotically stable%
\begin{equation}
x\left(  k+1\right)  =\left(  \left \{  A_{1},A_{2}\right \}  +\left \{
B_{1},B_{2}\right \}  \left \{  K_{1},K_{2}\right \}  \right)  x(k).
\label{closed}%
\end{equation}

\end{definition}

The following result was proven in \cite{zhou17arxiv}.

\begin{lemma}
\label{pp1}The complex-valued linear system (\ref{sys}) is stabilizable if and
only if%
\[
\mathrm{rank}\left[
\begin{array}
[c]{cccc}%
\lambda I_{n}-A_{1} & -A_{2}^{\#} & B_{1} & B_{2}^{\#}\\
-A_{2} & \lambda I_{n}-A_{1}^{\#} & B_{2} & B_{1}^{\#}%
\end{array}
\right]  =2n,\; \forall \lambda \in \{s:\left \vert s\right \vert \geq1\}.
\]

\end{lemma}

The following simple test for the stabilizability of the antilinear system
(\ref{antilinear}) was also recalled from \cite{zhou17arxiv}.

\begin{corollary}
\label{coro6}The antilinear system (\ref{antilinear}) is stabilizable if and
only if%
\begin{equation}
\mathrm{rank}\left[
\begin{array}
[c]{ccc}%
\lambda I_{n}-A_{2}A_{2}^{\#} & B_{2} & A_{2}B_{2}^{\#}%
\end{array}
\right]  =n,\; \forall \lambda \in \{s:\left \vert s\right \vert \geq1\},
\label{eq18a}%
\end{equation}
namely, the normal discrete-time linear system $(A_{2}A_{2}^{\#},[B_{2}%
,A_{2}B_{2}^{\#}])$ is stabilizable.
\end{corollary}

It follows that, for stabilization of the complex-valued linear system
(\ref{sys}), the full state feedback (\ref{eqfeedback}) is generally
necessary. However, for the discrete-time antilinear system (\ref{antilinear}%
), the well-used normal linear feedback%
\begin{equation}
u\left(  k\right)  =K_{1}x\left(  k\right)  , \label{normalfeedback}%
\end{equation}
is enough for stabilization under condition (\ref{eq18a}) \cite{zhou17arxiv}.

\subsection{Problem Formulation and Solution}

We study the linear quadratic regulation (LQR) problem for the complex-valued
linear system (\ref{sys}). Consider the real-valued quadratic index function%
\begin{equation}
J\left(  u\right)  =\sum \limits_{k=0}^{\infty}\left(  x^{\mathrm{H}%
}(k)Qx(k)+u^{\mathrm{H}}(k)Ru(k)\right)  , \label{eqj}%
\end{equation}
where $Q\in \mathbf{C}^{n\times n}$ and $R\in \mathbf{C}^{m\times m}$ are given
positive definite weighting matrices ($Q$ can be semi-positive definite,
however, we assume $Q>0$ for simplicity). The LQR problem refers to as finding
an optimal controller $u^{\ast}$ for system (\ref{sys}) such that $J\left(
u\right)  $ is minimized, denoted by $J_{\min}\left(  u^{\ast}\right)  $. The
LQR problem is said to be solvable if $J_{\min}\left(  u^{\ast}\right)
<\infty$ \cite{am07book,zhou17arxiv}.

The following result was proven in \cite{zhou17arxiv} regarding the existence
of a solution to the LQR problem.

\begin{lemma}
\label{lm2}The following statements are equivalent:

\begin{enumerate}
\item The LQR problem associated with system (\ref{sys}) has a solution.

\item The complex-valued linear system (\ref{sys}) is stabilizable.

\item There is a unique bimatrix $\left \{  P_{1},P_{2}\right \}  >0$ to the
following bimatrix Riccati equation%
\begin{align}
-\left \{  Q,0\right \}  =  &  \left \{  A_{1},A_{2}\right \}  ^{\mathrm{H}%
}\left \{  P_{1},P_{2}\right \}  \left \{  A_{1},A_{2}\right \}  -\left \{
P_{1},P_{2}\right \} \nonumber \\
&  -\left \{  A_{1},A_{2}\right \}  ^{\mathrm{H}}\left \{  P_{1},P_{2}\right \}
\left \{  B_{1},B_{2}\right \}  \left \{  S_{1},S_{2}\right \}  ^{-1}\left \{
B_{1},B_{2}\right \}  ^{\mathrm{H}}\left \{  P_{1},P_{2}\right \}  \left \{
A_{1},A_{2}\right \}  , \label{eqare}%
\end{align}
where $\left \{  S_{1},S_{2}\right \}  =\left \{  R,0\right \}  +\left \{
B_{1},B_{2}\right \}  ^{\mathrm{H}}\left \{  P_{1},P_{2}\right \}  \left \{
B_{1},B_{2}\right \}  .$
\end{enumerate}

Under one of the above conditions, the optimal control is the full state
feedback
\begin{equation}
u^{\ast}\left(  k\right)  =\left \{  K_{1}^{\ast},K_{2}^{\ast}\right \}
x\left(  k\right)  , \label{optcontrol}%
\end{equation}
where $\left \{  K_{1}^{\ast},K_{2}^{\ast}\right \}  $ is the optimal feedback
gain bimatrix determined by%
\begin{equation}
\left \{  K_{1}^{\ast},K_{2}^{\ast}\right \}  =-\left \{  S_{1},S_{2}\right \}
^{-1}\left \{  B_{1},B_{2}\right \}  ^{\mathrm{H}}\left \{  P_{1},P_{2}\right \}
\left \{  A_{1},A_{2}\right \}  , \label{eqgain}%
\end{equation}
the closed-loop system is asymptotically stable, and the minimal value of
$J\left(  u\right)  $ is given by
\begin{equation}
J_{\min}\left(  u^{\ast}\right)  =\operatorname{Re}\left(  x_{0}^{\mathrm{H}%
}\left \{  P_{1},P_{2}\right \}  x_{0}\right)  . \label{minj}%
\end{equation}

\end{lemma}

The above result is quite neat in the sense that the bimatrix Riccati equation
takes an analogous form as the usual Riccati matrix equation
\cite{kucera72kyb}.

\subsection{Iterative Solution to the Bimatrix Riccati Equation}

In this subsection, we provide an iterative method for solving the bimatrix
Riccati equation (\ref{eqare}). This method is not only useful for computing
solutions to (\ref{eqare}) but is also helpful in proving theoretical results
in the subsequent sections.

Motivated by the existing work for normal discrete-time Riccati equations
\cite{alks97natma}, we construct the following iteration associated with the
bimatrix Riccati equation (\ref{eqare}):%
\begin{align}
&  \left \{  P_{1}(k+1),P_{2}(k+1)\right \}  =\left \{  Q,0\right \}  +\left \{
A_{1},A_{2}\right \}  ^{\mathrm{H}}\left \{  P_{1}(k),P_{2}(k)\right \}  \left \{
A_{1},A_{2}\right \} \nonumber \\
&  -\left \{  A_{1},A_{2}\right \}  ^{\mathrm{H}}\left \{  P_{1}(k),P_{2}%
(k)\right \}  \left \{  B_{1},B_{2}\right \}  \left \{  S_{1}(k),S_{2}(k)\right \}
^{-1}\left \{  B_{1},B_{2}\right \}  ^{\mathrm{H}}\left \{  P_{1}(k),P_{2}%
(k)\right \}  \left \{  A_{1},A_{2}\right \}  , \label{eqit}%
\end{align}
where $\left \{  P_{1}(0),P_{2}(0)\right \}  =\left \{  Q,0\right \}  $ and
\[
\left \{  S_{1}(k),S_{2}(k)\right \}  =\left \{  R,0\right \}  +\left \{
B_{1},B_{2}\right \}  ^{\mathrm{H}}\left \{  P_{1}(k),P_{2}(k)\right \}  \left \{
B_{1},B_{2}\right \}  .
\]
For notation simplicity, we also denote
\begin{equation}
\left \{  R_{1},R_{2}\right \}  =\left \{  B_{1},B_{2}\right \}  \left \{
R,0\right \}  ^{-1}\left \{  B_{1},B_{2}\right \}  ^{\mathrm{H}}. \label{eqr}%
\end{equation}

\begin{theorem}
\label{th1}Assume that the complex-valued linear system (\ref{sys}) is
stabilizable and $\left \{  P_{1},P_{2}\right \}  $ is the unique positive
definite solution to (\ref{eqare}). Then, for any $k\geq0,$%
\begin{equation}
\left \{  Q,0\right \}  \leq \left \{  P_{1}(k),P_{2}(k)\right \}  \leq \left \{
P_{1}(k+1),P_{2}(k+1)\right \}  \leq \left \{  P_{1},P_{2}\right \}  .
\label{eq69}%
\end{equation}
Consequently, the limit of $\left \{  P_{1}(k),P_{2}(k)\right \}  $ as $k$
approaches infinity exists and%
\begin{equation}
\left \{  P_{1},P_{2}\right \}  =\lim_{k\rightarrow \infty}\left \{
P_{1}(k),P_{2}(k)\right \}  . \label{eq68}%
\end{equation}

\end{theorem}

By (\ref{eq69}) we can see that the iteration
(\ref{eqit}) can also be written as%
\begin{equation}
\left \{  P_{1}(k+1),P_{2}(k+1)\right \}  =\left \{  Q,0\right \}  +\left \{
A_{1},A_{2}\right \}  ^{\mathrm{H}}\left(  \left \{  P_{1}(k),P_{2}(k)\right \}
^{-1}+\left \{  R_{1},R_{2}\right \}  \right)  ^{-1}\left \{  A_{1}%
,A_{2}\right \}  , \label{eqitnew}%
\end{equation}
where $\left \{  P_{1}(0),P_{2}(0)\right \}  =\left \{  Q,0\right \}  $, and
$\left \{  R_{1},R_{2}\right \}  $ is given by (\ref{eqr}).

\section{\label{sec3}Optimal Control of Antilinear Systems}

In this section, we are interested in the antilinear system (\ref{antilinear}%
). Since it possesses a special structure, more specific results can be obtained.

\subsection{The Anti-Riccati Equation Based Approach}

We first present a so-called anti-Riccati equation based approach.

\begin{theorem}
\label{th2}Consider the antilinear system (\ref{antilinear}). Then the
following three statements are equivalent:

\begin{enumerate}
\item The LQR problem associated with system (\ref{antilinear}) has a solution.

\item The system (\ref{antilinear}) is stabilizable, namely, (\ref{eq18a}) is satisfied.

\item There is a unique positive definite solution $P_{A}>0$ to the so-called
anti-Riccati equation
\begin{equation}
-Q=A_{2}^{\mathrm{H}}P_{A}^{\#}A_{2}-A_{2}^{\mathrm{H}}P_{A}^{\#}B_{2}\left(
R+B_{2}^{\mathrm{H}}P_{A}^{\#}B_{2}\right)  ^{-1}B_{2}^{\mathrm{H}}P_{A}%
^{\#}A_{2}-P_{A}. \label{are2}%
\end{equation}

\end{enumerate}

In this case, the unique positive definite solutions to (\ref{eqare}) and
(\ref{are2}) are related with%
\begin{equation}
\left \{  P_{1},P_{2}\right \}  =\left \{  P_{A},0\right \}  . \label{eqp1p2pa}%
\end{equation}
Moreover, the optimal controller is the normal state feedback
(\ref{normalfeedback}) with $K_{1}=K_{1}^{\ast}$ defined by%
\begin{equation}
K_{1}^{\ast}=-\left(  R+B_{2}^{\mathrm{H}}P_{A}^{\#}B_{2}\right)  ^{-1}%
B_{2}^{\mathrm{H}}P_{A}^{\#}A_{2}, \label{eqk0}%
\end{equation}
the closed-loop system is asymptotically stable, and the optimal value of
$J\left(  u\right)  $ is
\begin{equation}
J_{\min}\left(  u\right)  =x_{0}^{\mathrm{H}}P_{A}x_{0}. \label{eqpja}%
\end{equation}

\end{theorem}

Theorem \ref{th2} improves some results in \cite{wqls16jfi} and
\cite{zhou17arxiv} where the existence of a positive definite solution to
(\ref{are2}) was not guaranteed. Moreover, we have relaxed controllability in
\cite{wqls16jfi} as stabilizability in this paper.

\begin{remark}
As a by-product of the proof of Theorem \ref{th2}, we can see
that the iteration%
\begin{equation}
P_{A}(k+1)=Q+A_{2}^{\mathrm{H}}\left(  P_{A}^{-\#}(k)+B_{2}R^{-1}%
B_{2}^{\mathrm{T}}\right)  ^{-1}A_{2}, \label{eqantiit}%
\end{equation}
with $P_{A}\left(  0\right)  =Q$, converges to the unique positive definite
solution to the anti-Riccati equation (\ref{are2}).
\end{remark}

Very recently, we have studied a class of nonlinear matrix equations in the
form of \cite{lzl14amc,zcl13amc}%
\begin{equation}
X+A^{\mathrm{H}}X^{-\#}A=I_{n}, \label{are3}%
\end{equation}
where $A\in \mathbf{C}^{n\times n}$ is known. Next we show how to link the
anti-Riccati equation (\ref{are2}) with this class of nonlinear matrix
equations. To this end, we define%
\begin{equation}
Q_{0}=Q^{-1}+\left(  A_{2}Q^{-1}A_{2}^{\mathrm{H}}\right)  ^{\#}+\left(
B_{2}R^{-1}B_{2}^{\mathrm{H}}\right)  ^{\#}>0. \label{eqq0}%
\end{equation}

\begin{proposition}
\label{lm3}If the anti-Riccati equation (\ref{are2}) has a positive definite
solution $P_{A},$ then the nonlinear matrix equation (\ref{are3}) with%
\begin{equation}
A=Q_{0}^{-\frac{\#}{2}}A_{2}Q^{-1}Q_{0}^{-\frac{1}{2}}, \label{eqa}%
\end{equation}
also has a positive definite solution $X$ such that%
\begin{equation}
X=Q_{0}^{-\frac{1}{2}}\left(  P_{A}^{-1}+\left(  A_{2}Q^{-1}A_{2}^{\mathrm{H}%
}\right)  ^{\#}+\left(  B_{2}R^{-1}B_{2}^{\mathrm{H}}\right)  ^{\#}\right)
Q_{0}^{-\frac{1}{2}}. \label{eqxx}%
\end{equation}
Moreover, if the antilinear system (\ref{antilinear}) is stabilizable, then
$X$ given by (\ref{eqxx}) is the maximal solution to (\ref{are3}).
\end{proposition}

By this proposition, when system (\ref{antilinear}) is stabilizable, the
unique positive solution to the anti-Riccati equation (\ref{are2}) can be
obtained by computing the maximal solution to the nonlinear matrix equation
(\ref{are3}) which has been carefully studied in \cite{lzl14amc} and
\cite{zcl13amc}. We finally remark that, as indicated by Proposition
\ref{lm3}, if system (\ref{antilinear}) is stabilizable and (\ref{are3}) has
any other positive definite solutions $X_{2},$ then we must have
\[
X_{2}\leq Q_{0}^{-\frac{1}{2}}\left(  \left(  B_{2}R^{-1}B_{2}^{\mathrm{H}%
}\right)  ^{\#}+\left(  A_{2}Q^{-1}A_{2}^{\mathrm{H}}\right)  ^{\#}\right)
Q_{0}^{-\frac{1}{2}}.
\]

\subsection{A Normal Riccati Equation Based Approach}

In this subsection, we establish a normal Riccati equation based approach to
the LQR problem for the antilinear system (\ref{antilinear}). For notation
simplicity, we denote
\begin{equation}
\left \{
\begin{array}
[c]{l}%
A_{N}=A_{2}^{\#}\left(  I_{n}-B_{2}\left(  R+B_{2}^{\mathrm{H}}Q^{\#}%
B_{2}\right)  ^{-1}B_{2}^{\mathrm{H}}Q^{\#}\right)  A_{2},\\
B_{N}=\left[
\begin{array}
[c]{cc}%
B_{2}^{\#} & A_{2}^{\#}B_{2}%
\end{array}
\right]  ,\\
Q_{N}=Q+A_{2}^{\mathrm{H}}\left(  Q^{-\#}+B_{2}R^{-1}B_{2}^{\mathrm{H}%
}\right)  ^{-1}A_{2},\\
R_{N}=\left[
\begin{array}
[c]{cc}%
R^{\#} & 0\\
0 & R+B_{2}^{\mathrm{H}}Q^{\#}B_{2}%
\end{array}
\right]  .
\end{array}
\right.  \label{eqanqr}%
\end{equation}

\begin{theorem}
\label{th4}Consider the antilinear system (\ref{antilinear}). Then the
following three statements are equivalent:

\begin{enumerate}
\item The LQR problem associated with system (\ref{antilinear}) has a solution.

\item The system (\ref{antilinear}) is stabilizable, namely, (\ref{eq18a}) is satisfied.

\item There is a unique positive definite solution $P_{N}$ to the normal
Riccati equation%
\begin{equation}
-Q_{N}=A_{N}^{\mathrm{H}}P_{N}A_{N}-P_{N}-A_{N}^{\mathrm{H}}P_{N}B_{N}\left(
R_{N}+B_{N}^{\mathrm{H}}P_{N}B_{N}\right)  ^{-1}B_{N}^{\mathrm{H}}P_{N}A_{N}.
\label{are5}%
\end{equation}

\end{enumerate}

In this case, the optimal controller is the normal state feedback
(\ref{normalfeedback}) with $K_{1}=K_{1}^{\ast}$ defined by%
\begin{equation}
K_{1}^{\ast}=-\left(  \left(  R+B_{2}^{\mathrm{H}}Q^{\#}B_{2}\right)
^{-1}B_{2}^{\mathrm{H}}Q^{\#}A_{2}+\left[
\begin{array}
[c]{cc}%
0 & I_{m}%
\end{array}
\right]  \left(  R_{N}+B_{N}^{\mathrm{H}}P_{N}B_{N}\right)  ^{-1}%
B_{N}^{\mathrm{H}}P_{N}A_{N}\right)  , \label{eqKN}%
\end{equation}
the closed-loop system is asymptotically stable and the minimal value of
$J\left(  u\right)  $ is%
\begin{equation}
J_{\min}\left(  u\right)  =x_{0}^{\mathrm{H}}P_{N}x_{0}. \label{eqjmin}%
\end{equation}

\end{theorem}

For system (\ref{antilinear}), Theorem \ref{th4} is better than Lemma
\ref{lm2} in the sense that the dimension of the Riccati equation in Theorem
\ref{th4} is half of that in Lemma \ref{lm2}, and Theorem \ref{th4} is better
than Theorem \ref{th2} in the sense that a normal Riccati equation is involved
in Theorem \ref{th4} while a non-standard Riccati equation is involved in
Theorem \ref{th2}.

Notice that we can construct the following iteration for the normal Riccati equation (\ref{are5}):
\begin{equation}
P_{N}\left(  k+1\right)  =Q_{N}+A_{N}^{\mathrm{H}}\left(  P_{N}^{-1}\left(
k\right)  +B_{N}R_{N}^{-1}B_{N}^{\mathrm{H}}\right)  ^{-1}A_{N},
\label{eqitnormal}%
\end{equation}
with $P_{N}\left(  0\right)  =Q_{N}.$  Then
\begin{equation}
\lim_{k\rightarrow \infty}P_{N}\left(  k\right)  =P_{N}. \label{eq58}%
\end{equation}

\subsection{Relationships Among Three Riccati Equations}

The following theorem links solutions to these three different Riccati
equations (\ref{eqare}), (\ref{are2}) and (\ref{are5}).

\begin{theorem}
\label{th3}Consider the antilinear system (\ref{antilinear}). Then the
following three statements are equivalent:

\begin{enumerate}
\item The bimatrix Riccati equation (\ref{eqare}) has a unique positive
definite solution $\{P_{1},P_{2}\}.$

\item The anti-Riccati equation (\ref{are2}) has a unique positive definite
solution $P_{A}.$

\item The normal Riccati equation (\ref{are5}) has a unique positive definite
solution $P_{N}.$
\end{enumerate}

Moreover, these solutions satisfy $P_{2}=0$ and%
\begin{equation}
P_{1}=P_{A}=P_{N}. \label{eqp1papn}%
\end{equation}

\end{theorem}

We can show that the
iteration (\ref{eqitnormal}) for the normal Riccati equation (\ref{are5})
converges faster than the iteration (\ref{eqantiit}) for the anti-Riccati
equation (\ref{are2}). Thus, from the computational point of view, the normal
Riccati equation (\ref{are5}) is recommended to use.

\begin{remark}
It follows from this theorem that the bimatrix Riccati equation based optimal
gain (\ref{eqgain}), the anti-Riccati equation based optimal gain
(\ref{eqk0}), and the normal Riccati equation based optimal gain (\ref{eqKN})
are equivalent.
\end{remark}

\begin{remark}
A simple proof for $P_{A}=P_{N}$ can be given as follows. By Theorems
\ref{th2} and \ref{th4}, the optimal problem is solvable with respectively the
minimal value $J_{\min}\left(  u\right)  =x_{0}^{\mathrm{H}}P_{A}x_{0}\ $and
$J_{\min}\left(  u\right)  =x_{0}^{\mathrm{H}}P_{N}x_{0}.$ As both $P_{A}$ and
$P_{N}$ are independent of $x_{0},$ we must have
\begin{equation}
x_{0}^{\mathrm{H}}P_{A}x_{0}=x_{0}^{\mathrm{H}}P_{N}x_{0},\; \forall x_{0}%
\in \mathbf{C}^{n}. \label{eq80}%
\end{equation}
Next we claim that, for positive definite matrices $P_{A}$ and $P_{N},$
$P_{A}=P_{N}$ if and only if (\ref{eq80}). Clearly, we need only to prove the
\textquotedblleft if\textquotedblright \ part. Denote $P_{A}=[a_{ij}]$ and
$P_{N}=[n_{ij}],i,j\in \mathbf{I}\left[  1,n\right]  .$ Letting $x_{0}=e_{i},$
where $e_{i}$ is the $i$th column of $I_{n},$ in (\ref{eq80}) gives
$a_{ii}=n_{ii},i\in \mathbf{I}\left[  1,n\right]  .$ Letting $x_{0}%
=[1,a+\mathrm{j}b,0,\ldots,0]^{\mathrm{H}},$ where $a\in \mathbf{R}%
,b\in \mathbf{R},$ in (\ref{eq80}) gives%
\[
a\operatorname{Re}\left(  a_{12}\right)  -b\operatorname{Im}\left(
a_{12}\right)  =a\operatorname{Re}\left(  n_{12}\right)  -b\operatorname{Im}%
\left(  n_{12}\right)  ,
\]
which, by respectively choosing $\left(  a=0,b\neq0\right)  $ and $\left(
b=0,a\neq0\right)  $, implies respectively $\operatorname{Im}\left(
a_{12}\right)  =\operatorname{Im}\left(  n_{12}\right)  $ and
$\operatorname{Re}\left(  a_{12}\right)  =\operatorname{Re}\left(
n_{12}\right)  ,$ namely, $a_{12}=n_{12}.$ Similarly, if we choose
$x_{0}=[1,0,a+\mathrm{j}b,0,\cdots,0]^{\mathrm{H}},$ we get $a_{13}=n_{13}.$
Repeating this process we finally have $P_{A}=P_{N}.$ However, the current
proof for Theorem \ref{th3} has its own value since it reveals the
relationship between the iteration (\ref{eqantiit}) for
the anti-Riccati equation (\ref{are2}) and the iteration (\ref{eqitnormal})
for the normal Riccati equation (\ref{are5}).
\end{remark}

\section{\label{sec4}Applications to Optimal Control of Time-Delay Systems}

\subsection{System and Problem Descriptions}

In this section, we consider the following discrete-time time-delay system
with only one step delay
\begin{equation}
\xi \left(  k+1\right)  =A_{0}\xi \left(  k\right)  +A_{\mathrm{d}}\xi \left(
k-1\right)  +Gv\left(  k\right)  ,\;k\geq0, \label{eqtds}%
\end{equation}
where $A_{0},A_{\mathrm{d}}\in \mathbf{R}^{n\times n}$ and $B\in \mathbf{R}%
^{n\times p}$ are known matrices, $\xi \in \mathbf{R}^{n}$ is the state vector,
and $v\in \mathbf{R}^{p}$ is the control vector. The initial condition is
$\xi \left(  0\right)  \in \mathbf{R}^{n}$ and $\xi \left(  -1\right)
\in \mathbf{R}^{n}.$ Without loss of generality, we assume that $p=2m,$ namely,
$p$ is an even number. Otherwise, we let $v=[v^{\mathrm{T}},w^{\mathrm{T}%
}]^{\mathrm{T}}$ and $G=[G,0]$£¬ where $w$ is any slack variable. Thus
we can let
\begin{equation}
G=\left[
\begin{array}
[c]{cc}%
G_{1} & G_{2}%
\end{array}
\right]  ,v\left(  k\right)  =\left[
\begin{array}
[c]{c}%
v_{1}\left(  k\right) \\
v_{2}\left(  k\right)
\end{array}
\right]  ,G_{i}\in \mathbf{R}^{n\times m},v_{i}\in \mathbf{R}^{m},i=1,2.
\label{eqg}%
\end{equation}

The problem to be solved is finding $v\left(  k\right)  $ for system
(\ref{eqtds}) such that the following quadratic index function is minimized%
\begin{equation}
J\left(  v\right)  =\sum \limits_{k=0}^{\infty}\left(  \xi^{\mathrm{T}}\left(
k\right)  Q_{0}\xi \left(  k\right)  +v^{\mathrm{T}}\left(  k\right)
R_{0}v\left(  k\right)  \right)  , \label{eqjdelay}%
\end{equation}
where $Q_{0}\in \mathbf{R}^{n\times n}$ and $R_{0}\in \mathbf{R}^{m\times m}$
are given positive definite matrices.

\begin{remark}
We explain that we can assume without loss of generality that $R_{0}$ is a
block diagonal matrix. Denote%
\[
R_{0}=\left[
\begin{array}
[c]{cc}%
R_{01} & R_{02}\\
R_{02}^{\mathrm{T}} & R_{03}%
\end{array}
\right]  ,\;R_{0i}\in \mathbf{R}^{m\times m},i=1,2,3.
\]
Since $R_{0}>0,$ by the Schur complement, we have $R_{03}-R_{02}^{\mathrm{T}%
}R_{01}^{-1}R_{02}>0.$ Thus we can denote%
\[
L_{0}=\left[
\begin{array}
[c]{cc}%
I_{m} & -R_{01}^{-1}R_{02}\left(  R_{03}-R_{02}^{\mathrm{T}}R_{01}^{-1}%
R_{02}\right)  ^{-\frac{1}{2}}R_{01}^{\frac{1}{2}}\\
0 & \left(  R_{03}-R_{02}^{\mathrm{T}}R_{01}^{-1}R_{02}\right)  ^{-\frac{1}%
{2}}R_{01}^{\frac{1}{2}}%
\end{array}
\right]  .
\]
Direct computation gives%
\[
L_{0}^{\mathrm{T}}R_{0}L_{0}=\left[
\begin{array}
[c]{cc}%
R_{01} & \\
& R_{01}%
\end{array}
\right]  >0.
\]
Then, by the input transformation $\hat{v}=L_{0}v,$ the time-delay system
(\ref{eqtds}) can be written as%
\[
\xi \left(  k+1\right)  =A_{0}\xi \left(  k\right)  +A_{\mathrm{d}}\xi \left(
k-1\right)  +\hat{G}\hat{v}\left(  k\right)  ,
\]
where $\hat{G}=GL_{0},$ and the quadratic index function (\ref{eqjdelay})
becomes%
\begin{align*}
J\left(  v\right)   &  =\sum \limits_{k=0}^{\infty}\left(  \xi^{\mathrm{T}%
}\left(  k\right)  Q_{0}\xi \left(  k\right)  +\hat{v}^{\mathrm{T}}\left(
k\right)  L_{0}^{\mathrm{T}}R_{0}L_{0}\hat{v}\left(  k\right)  \right) \\
&  =\sum \limits_{k=0}^{\infty}\left(  \xi^{\mathrm{T}}\left(  k\right)
Q_{0}\xi \left(  k\right)  +\hat{v}^{\mathrm{T}}\left(  k\right)  \left[
\begin{array}
[c]{cc}%
R_{01} & \\
& R_{01}%
\end{array}
\right]  \hat{v}\left(  k\right)  \right)  .
\end{align*}
Thus, without loss of generality, we can assume that%
\begin{equation}
R_{0}=\left[
\begin{array}
[c]{cc}%
R & \\
& R
\end{array}
\right]  ,\;0<R\in \mathbf{R}^{m\times m}. \label{eqr0}%
\end{equation}

\end{remark}

Therefore we assume hereafter that $R_{0}$ takes the special form (\ref{eqr0}).

\begin{proposition}
The time-delay system (\ref{eqtds}) can be equivalently written as (\ref{sys})
where
\begin{equation}
\left \{
\begin{array}
[c]{l}%
x\left(  k\right)  =\xi \left(  k\right)  +\mathrm{j}\xi \left(  k-1\right)  ,\\
u\left(  k\right)  =v_{1}\left(  k\right)  +\mathrm{j}v_{2}\left(  k\right)
,\;k\geq0,
\end{array}
\right.  \label{eqxu}%
\end{equation}
and
\begin{equation}
\left \{
\begin{array}
[c]{cc}%
A_{1}=\frac{1}{2}A_{0}+\frac{\mathrm{j}}{2}\left(  I_{n}-A_{\mathrm{d}%
}\right)  , & B_{1}=\frac{1}{2}G_{1}-\frac{\mathrm{j}}{2}G_{2},\\
A_{2}=\frac{1}{2}A_{0}-\frac{\mathrm{j}}{2}\left(  I_{n}+A_{\mathrm{d}%
}\right)  , & B_{2}=\frac{1}{2}G_{1}-\frac{\mathrm{j}}{2}G_{2}.
\end{array}
\right.  \label{eqab}%
\end{equation}
Moreover, if $R_{0}$ takes the form (\ref{eqr0}), the quadratic index function
(\ref{eqjdelay}) can be written as
\begin{equation}
J_{1}\left(  u\right)  =\sum \limits_{k=0}^{\infty}\left(  x^{\mathrm{H}%
}\left(  k\right)  Qx\left(  k\right)  +u^{\mathrm{H}}\left(  k\right)
Ru\left(  k\right)  \right)  -\xi^{\mathrm{T}}\left(  -1\right)  Q\xi \left(
-1\right)  , \label{eqjdelay1}%
\end{equation}
where $Q=\frac{1}{2}Q_{0}>0.$
\end{proposition}

Since the last term $\xi^{\mathrm{T}}\left(  -1\right)  Q\xi \left(  -1\right)
$ in (\ref{eqjdelay1}) dependents on only the initial condition, $J_{1}\left(
u\right)  $ is minimized if and only if%
\begin{equation}
J_{2}\left(  u\right)  =\sum \limits_{k=0}^{\infty}\left(  x^{\mathrm{H}%
}\left(  k\right)  Qx\left(  k\right)  +u^{\mathrm{H}}\left(  k\right)
Ru\left(  k\right)  \right)  , \label{eqj2}%
\end{equation}
is minimized. Hence, the linear optimal control problem for the time-delay
system (\ref{eqtds}) has been transformed equivalently to the linear quadratic
optimal control problem for the complex-valued linear system (\ref{sys}) with
the quadratic index function (\ref{eqj2}). According to results in Section
\ref{sec2}, the solution to this problem has been completely characterized by
Lemma \ref{lm2}. Thus the optimal control is $u\left(  k\right)  =K_{1}^{\ast
}x\left(  k\right)  +\left(  K_{2}^{\ast}\right)  ^{\#}x^{\#},$ which, by
separating real and imaginary parts, is equivalent to \cite{zhou17arxiv}
\[
v\left(  k\right)  =\left[
\begin{array}
[c]{c}%
v_{1}\left(  k\right) \\
v_{2}\left(  k\right)
\end{array}
\right]  =\left[
\begin{array}
[c]{cc}%
\mathrm{Re}\left(  K_{1}^{\ast}+K_{2}^{\ast}\right)  & -\mathrm{Im}\left(
K_{1}^{\ast}+K_{2}^{\ast}\right) \\
\mathrm{Im}\left(  K_{1}^{\ast}-K_{2}^{\ast}\right)  & \mathrm{Re}\left(
K_{1}^{\ast}-K_{2}^{\ast}\right)
\end{array}
\right]  \left[
\begin{array}
[c]{c}%
\xi \left(  k\right) \\
\xi \left(  k-1\right)
\end{array}
\right]  ,
\]
which is physically implementable \cite{zhou17arxiv}.

\subsection{An Illustrative Example}

In this subsection we use the linearized F-16 aircraft model studied
previously in \cite{ss85jgcd} and \cite{lz16ietcta} to illustrate the obtained
results. The continuous-time model is shown as follows%
\begin{equation}
\dot{\xi}\left(  t\right)  =\mathcal{A}\xi \left(  t\right)  +\mathcal{A}%
_{\mathrm{d}}\xi \left(  t-\tau \right)  +\mathcal{G}v\left(  t\right)  ,
\label{a1}%
\end{equation}
in which we have assumed that there is a state delay $\tau=0.1$ in the
elevator deflection which is the fourth element of $x(t)$ \cite{lz16ietcta}.
The coefficient matrices are then given by \cite{lz16ietcta}%
\begin{align*}
\mathcal{A}  &  =\left[
\begin{array}
[c]{ccccc}%
0 & 1.0 & 0 & 0 & 0\\
0 & -0.8694 & 43.223 & -17.251 & -1.5766\\
0 & 0.9934 & -1.3411 & -0.1690 & -0.2518\\
0 & 0 & 0 & 0 & 0\\
0 & 0 & 0 & 0 & -20.0
\end{array}
\right]  ,\\
\mathcal{A}_{\mathrm{d}}  &  =\left[
\begin{array}
[c]{ccccc}%
0 & 0 & 0 & 0 & 0\\
0 & 0 & 0 & 0 & 0\\
0 & 0 & 0 & 0 & 0\\
0 & 0 & 0 & -20.0 & 0\\
0 & 0 & 0 & 0 & 0
\end{array}
\right]  ,\mathcal{G}=\left[
\begin{array}
[c]{cc}%
0 & 0\\
0 & 0\\
0 & 0\\
20.0 & 0\\
0 & 20.0
\end{array}
\right]  .
\end{align*}
By taking the sampling period as $T=0.1\mathrm{s}$, the continuous-time
time-delay system (\ref{a1}) can be discretized as (\ref{eqtds}) where%
\begin{align*}
A_{0}  &  =\left[
\begin{array}
[c]{ccccc}%
1.0000 & 0.1025 & 0.2080 & -0.0879 & -0.0057\\
0 & 1.1175 & 4.1534 & -1.8042 & -0.1010\\
0 & 0.0955 & 1.0722 & -0.0994 & -0.0153\\
0 & 0 & 0 & 1.0000 & 0\\
0 & 0 & 0 & 0 & 0.1353
\end{array}
\right]  ,\; \\
A_{\mathrm{d}}  &  =\left[
\begin{array}
[c]{ccccc}%
0 & 0 & 0 & 0.0594 & 0\\
0 & 0 & 0 & -1.8165 & 0\\
0 & 0 & 0 & 0.0434 & 0\\
0 & 0 & 0 & -2.0000 & 0\\
0 & 0 & 0 & 0 & 0
\end{array}
\right]  ,G=\left[
\begin{array}
[c]{cc}%
-0.0581 & -0.0040\\
-1.7586 & -0.1131\\
-0.0720 & -0.0175\\
2.0000 & 0\\
0 & 0.8647
\end{array}
\right]  .
\end{align*}

We now consider the corresponding linear optimal control problem
(\ref{eqjdelay}) with $Q_{0}=2I_{5}$ and $R=1.$ Thus $Q=I_{5}.$ Let $\left \{
P_{1}(k),P_{2}(k)\right \}  $ be computed according to the iteration
(\ref{eqitnew}). Denote $e\left(  k\right)  =\ln \left \Vert \left \{
E_{1}\left(  k\right)  ,E_{2}\left(  k\right)  \right \}  \right \Vert $ where%
\[
\left \{  E_{1}\left(  k\right)  ,E_{2}\left(  k\right)  \right \}  =\left \{
A_{1},A_{2}\right \}  ^{\mathrm{H}}\left(  \left \{  P_{1}(k),P_{2}(k)\right \}
^{-1}+\left \{  R_{1},R_{2}\right \}  \right)  ^{-1}\left \{  A_{1}%
,A_{2}\right \}  +\left \{  Q,0\right \}  -\left \{  P_{1}(k),P_{2}(k)\right \}  .
\]
Numerical computation indicates that $e(k)$ converges to zero in a rather fast
speed. For $k=140,$ we obtain%
\begin{align*}
P_{1}\left(  k\right)   &  =%
\begin{bmatrix}%
\begin{smallmatrix}
12.3464+0.0000\mathrm{j} & 2.3671+0.0000\mathrm{j} & 8.6342+0.0000\mathrm{j} &
-2.0194+3.2919\mathrm{j} & -0.1762+0.0000\mathrm{j}\\
2.3671+0.0000\mathrm{j} & 3.7958+0.0000\mathrm{j} & 10.5820+0.0000\mathrm{j} &
-3.6990+5.5020\mathrm{j} & -0.2025+0.0000\mathrm{j}\\
8.6342+0.0000\mathrm{j} & 10.5820+0.0000\mathrm{j} & 56.7625+0.0000\mathrm{j}
& -18.9340+23.6474\mathrm{j} & -0.9975+0.0000\mathrm{j}\\
-2.0194-3.2919\mathrm{j} & -3.6990-5.5020\mathrm{j} &
-18.9340-23.6474\mathrm{j} & 28.0046+0.0000\mathrm{j} &
0.3487+0.4725\mathrm{j}\\
-0.1762+0.0000\mathrm{j} & -0.2025+0.0000\mathrm{j} & -0.9975+0.0000\mathrm{j}
& 0.3487-0.4725\mathrm{j} & 1.5259+0.0000\mathrm{j}%
\end{smallmatrix}
\end{bmatrix}
,\\
P_{2}\left(  k\right)   &  =%
\begin{bmatrix}%
\begin{smallmatrix}
11.3464+0.0000\mathrm{j} & 2.3671+0.0000\mathrm{j} & 8.6342+0.0000\mathrm{j} &
-2.0194+3.2919\mathrm{j} & -0.1762+0.0000\mathrm{j}\\
2.3671+0.0000\mathrm{j} & 2.7958+0.0000\mathrm{j} & 10.5820+0.0000\mathrm{j} &
-3.6990+5.5020\mathrm{j} & -0.2025+0.0000\mathrm{j}\\
8.6342+0.0000\mathrm{j} & 10.5820+0.0000\mathrm{j} & 55.7625+0.0000\mathrm{j}
& -18.9340+23.6474\mathrm{j} & -0.9975+0.0000\mathrm{j}\\
-2.0194+3.2919\mathrm{j} & -3.6990+5.5020\mathrm{j} &
-18.9340+23.6474\mathrm{j} & -3.5897-17.1002\mathrm{j} &
0.3487-0.4725\mathrm{j}\\
-0.1762+0.0000\mathrm{j} & -0.2025+0.0000\mathrm{j} & -0.9975+0.0000\mathrm{j}
& 0.3487-0.4725\mathrm{j} & 0.5259+0.0000\mathrm{j}%
\end{smallmatrix}
\end{bmatrix}
.
\end{align*}
It follows that $P_{1}\left(  k\right)  =P_{1}^{\mathrm{H}}\left(  k\right)  $
and $P_{2}\left(  k\right)  =P_{2}^{\mathrm{T}}\left(  k\right)  .$
Consequently, the optimal feedback gain $\{K_{1}^{\ast},K_{2}^{\ast}\}$ can be
computed according to (\ref{eqgain}) as%
\begin{align*}
K_{1}^{\ast}  &  =\left[
\begin{array}
[c]{ccccc}%
0.0463+0.0962\mathrm{j} & 0.1140+0.1205\mathrm{j} & 0.6384+0.6279\mathrm{j} &
-0.7529-0.4637\mathrm{j} & -0.0112-0.0584\mathrm{j}%
\end{array}
\right]  ,\\
K_{2}^{\ast}  &  =\left[
\begin{array}
[c]{ccccc}%
0.0463-0.0962\mathrm{j} & 0.1140-0.1205\mathrm{j} & 0.6384-0.6279\mathrm{j} &
-0.2122-0.0036\mathrm{j} & -0.0112+0.0584\mathrm{j}%
\end{array}
\right]  .
\end{align*}
Finally, with the initial condition%
\[
\left[
\begin{array}
[c]{c}%
\xi^{\mathrm{T}}\left(  0\right) \\
\xi^{\mathrm{T}}\left(  -1\right)
\end{array}
\right]  =\left[
\begin{array}
[c]{ccccc}%
4 & 1 & -8 & -6 & 9\\
4 & 4 & 8 & -6 & 10
\end{array}
\right] ,
\]
by simulation we can observe the asymptotic stability of the
closed-loop system.

\section{\label{sec5}Conclusion}

This paper has studied linear optimal control (linear quadratic regulation,
LQR for short) of discrete-time complex-valued linear systems. Firstly, an
iterative algorithm was proposed to solve the associated bimatrix Riccati
equation introduced in our early study. The convergence of the algorithm was
proven. Then the LQR problem for the antilinear system, which is a special
case of the complex-valued linear system, was carefully studied and three
different solutions were obtained, namely, bimatrix Riccati equation based
solution, anti-Riccati equation based solution, and normal Riccati equation
based solution. Relationships among these three different solutions are
revealed. The bimatrix Riccati equation based approach was then used to solve
the LQR problem of linear time-delay systems with one step state delay and an
illustrative example demonstrated the effectiveness of the proposed approach.


\bigskip

\begin{thebibliography}{99}                                                                                               %

\bibitem {am07book}Anderson B D O, Moore J B. \emph{Optimal Control: Linear
Quadratic Methods}. Dover Publications, 2007.

\bibitem {alks97natma}Assimakis N D, Lainiotis D G, Katsikas S K, Sanida F L.
A survey of recursive algorithms for the solution of the discrete time Riccati
equation. \textit{Nonlinear Analysis: Theory, Methods \& Applications}, 1997,
30(4): 2409-2420.

\bibitem {dl71tac}Dorato P, Levis A. Optimal linear regulators: The
discrete-time case. \emph{IEEE Transactions on Automatic Control}, 1971,
16(6): 613-620.

\bibitem {kailath80book}Kailath T. \textit{Linear Systems}, Englewood Cliffs,
NJ: Prentice-Hall, 1980.

\bibitem {kucera72kyb}Ku\v{c}era V. The discrete Riccati equation of optimal
control, \textit{Kybernetika}, 1972, 8(5): (430)-447.

\bibitem {lzl14amc}Li Z Y, Zhou B, Lam J. Towards positive definite solutions
of a class of nonlinear matrix equations, \textit{Applied Mathematics and
Computation}, 2014, 237: 546-559.

\bibitem {lz16ietcta}Liu Q, Zhou B. Delay compensation of discrete-time linear
systems by nested prediction. \emph{IET Control Theory \& Applications}, 2016,
10(15): 1824-1834.

\bibitem {rugh96book}Rugh W J. \textit{Linear System Theory}, Upper Saddle
River, NJ: prentice hall, 1996.

\bibitem {ss85jgcd}Sobel K M, Shapiro E Y. A design methodology for pitch
pointing flight control systems. \textit{Journal of Guidance, Control, and
Dynamics}, 1985, 8(2): 181-187.

\bibitem {woodbury49book}Woodbury MA. \textit{The Stability of Out-input
Matrices}. University of Chicago Press, Chicago, p. 93, 1949.

\bibitem {wdl13aucc}Wu A G, Duan G R, Liu W, Sreeram V. Controllability and
stability of discrete-time antilinear systems. \textit{2013 3rd Australian
Control Conference (AUCC)}, 2013: 403-408.

\bibitem {wqls16jfi}Wu A G, Qian Y Y, Liu W, Sreeram V. Linear quadratic
regulation for discrete-time antilinear systems: An anti-Riccati matrix
equation approach. \textit{Journal of the Franklin Institute}, 2016, 353(5): 1041-1060.

\bibitem {wzls15iet}Wu A G, Zhang Y, Liu W, Sreeram V. State response for
continuous-time antilinear systems. \textit{IET Control Theory \&
Applications}, 2015, 9(8): 1238-1244.

\bibitem {zhang17auto}Zhang G. Dynamical analysis of quantum linear systems
driven by multi-channel multi-photon states, \textit{Automatica}, 2017, 83: 186-198.

\bibitem {zhou17arxiv}Zhou B. Analysis and design of complex-valued linear
systems. arXiv preprint arXiv:1708.05120, 2017.

\bibitem {zhou17arxiv2}Zhou B. Solutions to linear bimatrix equations with
applications to pole assignment of complex-valued linear systems. arXiv
preprint arXiv:1708.07947, 2017.

\bibitem {zcl13amc}Zhou B, Cai G B, Lam J. Positive definite solutions of the
nonlinear matrix equation $X+A^{\mathrm{H}}\overline{X}^{-1}A=I$,
\textit{Applied Mathematics and Computation}, 2013, 219(14): 7377-7391.
\end{thebibliography}
\end{document}